\documentclass[12pt,oneside]{amsart}
\usepackage{amsmath,amsfonts,latexsym,graphicx,amssymb,url}
\usepackage{amsmath,txfonts,pifont}
\newcommand{\diam}{\text{diam}}
\newcommand{\osc}{\text{\rm osc}}
\newcommand{\supp}{\text{{\rm supp }}}
\newcommand{\trace}{\text{\rm trace\,}}
\newcommand{\R}{{\mathbb R}} 
\newcommand{\C}{{\mathbb C}}
\newcommand{\N}{{\mathbb N}}
\def \H {{\mathbb H}}
\def \e {\varepsilon}
\def \O {\Omega}
\def \l {\lambda}
\def \p {\partial}

\def \s {\sigma}

\theoremstyle{plain}
\newtheorem{theorem}{Theorem}[section]
\newtheorem{corollary}[theorem]{Corollary}
\newtheorem{lemma}[theorem]{Lemma}
\newtheorem{proposition}[theorem]{Proposition}
\newtheorem{definition}[theorem]{Definition}
\newtheorem{remark}[theorem]{Remark}

\begin{document}

\setcounter{equation}{0}










\title[Maximum and comparison principles on the Heisenberg group]{Maximum and comparison principles
for convex functions on the Heisenberg group}
\author[C. E. Guti\'errez and A. Montanari]{Cristian E.
Guti\'errez\\
\\
and\\ \\
Annamaria Montanari\\
\\
}
\thanks{November 13, 2002.\\The results in this paper and the ideas of their proofs have been presented
in the following talks: Analysis Seminar, Temple U., October 2002;
Fabes--Chiarenza Lectures at Siracusa, December 2002; Pan-American
Conference, Santiago de Chile, January 2003, Analysis Seminar, U. of Bologna, March 2003,
and Analysis Seminar,
U. Texas at Austin, March 2003.\\The first author was partially
supported by NSF grant DMS--0070648, and thanks the University of Bologna for its support
and hospitality in his several visits to carry out this project.
The second author thanks Temple University for the
hospitality during her visit during March--April 2001 when the key Theorem 3.1 was proved.}
\address{Department of Mathematics\\
    Temple University\\
    Philadelphia,
     PA 19122}

\email{gutierrez@math.temple.edu}

\address{Dipartimento di Matem\`atica\\
    Universit\`a di Bologna\\
    Piazza Porta San Donato 5\\
   Bologna, 40127,  Italy}

\email{montanar@dm.unibo.it}

\maketitle

\setcounter{equation}{0}
\section{Introduction}
The purpose in this paper is to establish pointwise estimates for
a class of convex functions on the Heisenberg group. An integral
estimate for classical convex functions in terms of the
Monge--Amp\`ere operator $\det D^2u$ was proved by Aleksandrov,
see \cite[Theorem 1.4.2]{Gut:book}. Such estimate is of great
importance in the theory of weak solutions for the Monge--Amp\`ere
equation, and its proof revolves around the geometric features of
the notion of normal mapping or subdifferential in $\R^n$
\cite[Definition 1.1.1]{Gut:book} which yield in addition the
useful comparison principle for Monge-Amp\`ere measures,
\cite[Theorem 1.4.6]{Gut:book}.

On the Heisenberg group, and more generally in Carnot groups,
several notions of convexity have been introduced and compared in
\cite{Danielli-Garofalo-Nhieu:notionsofconvexity} and
\cite{Lu-Manfredi-Stroffolini:notionsofconvexityinheisenberg}. The
notion of convex function we use in this paper is given in
Definition \ref{def:convexity}, and a natural question is if
similar comparison and maximum principles hold in this setting. A
reason for this question is that those estimates would be useful
in the study of solutions for nondivergence equations of the form
$a_{ij}X_i\,X_j$ where $a_{ij}$ is a uniformly elliptic measurable
matrix and $X_i$ are the Heisenberg vector fields. The difficulty
for this study is the doubtful existence of a notion of normal
mapping in $\H^n$ suitable to establish maximum and comparison
principles.

In this paper we address this question and follow a route
different from the one described above for convex functions, and
in particular, we do not use any notion of normal mapping. This
approach was recently used by Trudinger and Wang to study Hessian
equations \cite{Trudinger-Wang:hessianmeasuresI}. Our integral
estimates are in terms of the following Monge--Amp\`ere type
operator: $\det \mathcal H(u) + 12\,(u_t)^2$, see Definition
\ref{def:Hconvex}. We first establish by means of integration by
parts a comparison principle for smooth functions, Theorem
\ref{thm:comparisonprinciple}, and then extend this principle to
"cones" Theorem \ref{thm:comparisonprincipleperforated}. This
together with the geometry in $\H^n$ leads by iteration to the
maximum principle Theorem \ref{ABP}. We next estimate the
oscillation of $\mathcal H$--convex functions Proposition
\ref{pro:osc} that permits to extend our definition of
Monge--Amp\`ere measure to continuous $\mathcal H$--convex
functions and obtain a general comparison principle Theorem
\ref{thm:comparisonpcpleformeasures}.

The paper is organized as follows. Section \ref{sec:preliminaries}
contains preliminaries about $\H^n$ and the definitions of
$\mathcal H$--convexity. In Section \ref{sec:comparisonprinciple}
we prove the comparison principle for $C^2$ functions. Section
\ref{sec:weakmaxprinciple} contains the proof that "cones"
agreeing with $\mathcal H$--convex functions $u$ on the boundary
are above $u$ inside, and the comparison principle for cones
Theorem \ref{thm:comparisonprincipleperforated}. In Section
\ref{sec:maximumprinciple} we prove a maximum principle similar to
Aleksandrov's estimate aforementioned. Finally, Section
\ref{sec:Hmeasures} contains the oscillation estimates and the
construction of the analogue of Monge--Amp\`ere measures for
$\mathcal H$--convex functions.

\setcounter{equation}{0}
\section{Preliminaries and $\mathcal
H$--convexity}\label{sec:preliminaries}









Let $u=u(x,y,t)$; $z=(x,y,t)$, and $X= \partial_x + 2y\,
\partial_t$, $Y=\partial_y - 2x\,
\partial_t$.
We have $[X,Y]=XY-YX=- 4\partial_t.$ If $ \xi_0=(x_0,y_0,t_0)$ and
$\xi=(x,y,t)$, then the non--commutative multiplication law in
$\H^1$
is given by
\[
\xi_0\circ \xi=(x_0 + x,y_0 +y, t_0+t+2(xy_0-yx_0)),
\]
and we have $\xi^{-1}=-\xi$, and $(\xi_0\circ
\xi)^{-1}=\xi^{-1}\circ  \xi_0^{-1}$. The gauge in $\H^1$ is
\[
\rho(\xi)=\left((x^2+y^2)^2+t^2\right)^{1/4},
\]
and the distance
\[
d(\xi,\xi_0)=\rho(\xi_0^{-1}\circ \xi).
\]
We have
\begin{equation}\label{triangle}
d(\xi,\xi_0)\leq d(\xi,\zeta)+d(\zeta,\xi_0)
\end{equation}
for every $\xi,\xi_0,\zeta\in \H^1.$ Given $\l>0$ we consider the
dilations
\[
\delta_{\lambda}(\xi)=(\l x, \l y, \l^2 t).
\]
Then
\[
d(\delta_\l\xi,\delta_\l \xi_0)=\l \,d(\xi,\xi_0).
\]
For more details about $\H^n$ see \cite[Chapters XII and
XIII]{S:book}.

\subsection{$\mathcal
H$--convexity}

Let $\xi_0=(x_o,y_o,t_o), \zeta=(x,y,t)$ and
\[
g(\zeta)=f(\xi_o\circ \zeta).
\]
We have
\[
\partial_xg(0)= Xf(\xi_o),
\quad
\partial_yg(0)= Yf(\xi_o),
\quad
\partial_tg(0)= \partial_tf(\xi_o),
\]
and
\[
\partial_{xx}g(0)=(X^2f)(\xi_0),
\quad
\partial_{xy}g(0)=(YXf)(\xi_0) - 2 \,\partial_tf(\xi_o),
\quad
\partial_{xt}g(0)=\partial_{tx}f (\xi_0) + 2 \,y_0\, \partial_{tt}f(\xi_o),
\]
\[
\partial_{yx}g(0)=(XYf)(\xi_0) + 2 \,\partial_tf(\xi_o),
\quad
\partial_{yy}g(0)=(Y^2f)(\xi_0),
\quad
\partial_{yt}g(0)=\partial_{ty}f (\xi_0) - 2 \,x_0\, \partial_{tt}f(\xi_o),
\]
\[
\partial_{tx}g(0)=\partial_{tx}f (\xi_0) + 2 \,y_0\, \partial_{tt}f(\xi_o),
\quad
\partial_{ty}g(0)=\partial_{ty}f (\xi_0) - 2 \,x_0\, \partial_{tt}f(\xi_o),
\quad
\partial_{tt}g(0)=\partial_{tt}f (\xi_0).
\]
Let
\[
A=
\left[\begin{matrix}
(X^2f)(\xi_0) & (YXf)(\xi_0) - 2 \,\partial_tf(\xi_o) &
\partial_{tx}f (\xi_0) + 2 \,y_0\, \partial_{tt}f(\xi_o)\\
(XYf)(\xi_0) + 2 \,\partial_tf(\xi_o) &
(Y^2f)(\xi_0) & \partial_{ty}f (\xi_0) - 2 \,x_0\,
\partial_{tt}f(\xi_o)\\
\partial_{tx}f (\xi_0) + 2 \,y_0\, \partial_{tt}f(\xi_o)
&
\partial_{ty}f (\xi_0) - 2 \,x_0\, \partial_{tt}f(\xi_o)
&
\partial_{tt}f (\xi_0)
\end{matrix}\right].
\]
Then the Taylor polynomial of order two of $g$ is
\begin{align*}
&f(\xi_o)
+
(Xf(\xi_o),Yf(y_0), \partial_tf(\xi_0))\cdot \zeta
+
\frac12
\langle A \zeta,\zeta \rangle\\
&=
f(\xi_o)
+
(Xf(\xi_o),Yf(y_0))\cdot (x,y)
+
(X^2f)\, x^2 + (XYf + YX f)\,xy + (Y^2f) \, y^2\\
&\qquad +
t\{
f_t(\xi_0)
+
2f_{tx} \, x + 4 y_o f_{tt} x +
(f_{tx}+f_{ty})y
- 4 x_o f_{tt} y\}.
\end{align*}
That is, if $(x,y,t)\in \Pi_0$ then $t=0$ and so
on this plane we have
\begin{align*}
g(\zeta)
&=
f(\xi_o)
+
(Xf(\xi_o),Yf(y_0))\cdot (x,y)\\
&\qquad +
(X^2f)\, x^2 + (XYf + YX f)\,xy + (Y^2f) \, y^2 +o(x^2+y^2).
\end{align*}

Set $B_R(\xi_0)=\{\xi \in \R^3: d(\xi,\xi_0)<R\}$. Given
$\xi_0=(x_0,y_0,t_0)\in R^3$ let
\[
\Pi_{\xi_0}=\{(x,y,t): t-t_0-2(xy_0-yx_0)=0 \}.
\]
That is, $\Pi_{\xi_0}$ is the plane generated by the vectors
$(1,0,2y_0)$, $(0,1,-2x_0)$ and passing through the point $\xi_0$.
Notice that if $h\in \H^1$, then
\begin{equation}\label{eq:invarianceplanes}
\xi\in \Pi_{\xi_0} \text{ if and only if } h\circ \xi\in
\Pi_{h\circ \xi_0}.
\end{equation}

Given $c\in \C$ and $u\in C^2(\Omega)$, let
\[
\mathcal H_c(u)= \left[\begin{matrix} X^2u & XYu + c u_t \\
YXu - c u_t  & Y^2u
\end{matrix}\right]
\]
and
\begin{equation}\label{eq:defofHc}
H_c(u)= \det
\left[\begin{matrix} X^2u & XYu + c u_t \\
YXu - c u_t  & Y^2u
\end{matrix}\right].
\end{equation}
\begin{definition}\label{def:Hconvex}
The function $u\in C^2(\Omega)$ is $\mathcal H$--convex in
$\Omega$ if the symmetric matrix
\[
\mathcal H(u)=\mathcal H_2(u)=\left[\begin{matrix} X^2u & (XYu + YXu)/2 \\
(XYu + YXu)/2  & Y^2u
\end{matrix}\right]
\]
is positive semidefinite in $\Omega$.
\end{definition}
Notice that the matrix $\mathcal H_c(u)$ is symmetric if and only
if $c=2.$ Also, if $\langle \mathcal H_c(u) \xi,\xi \rangle\geq 0$
for all $\xi\in \R^2$ and for some $c$, then this quadratic form
is nonnegative for all values of $c\in \R$.

We extend the definition of ${\mathcal H}$--convexity to
continuous functions.

\begin{definition}\label{def:convexity}
The function $u\in C(\Omega)$ is ${\mathcal H}$--convex in
$\Omega$ if there exists a sequence $u_k\in C^2(\Omega)$ of
${\mathcal H}$--convex functions in $\Omega$ such that $u_k\to u$
uniformly on compact subsets of $\Omega$.
\end{definition}

The following proposition yields equivalent definitions of
${\mathcal H}$--convexity, see \cite[Theorem 5.11]{Danielli-Garofalo-Nhieu:notionsofconvexity}
for the proof.

\begin{proposition}\label{prop:H-conveximpliesconvexonlines}
Let $u\in C(\Omega)$ with $\Omega\subset \R^3$ open\footnote{We
assume that if $\xi,\xi_0\in \Omega$, then $\xi_0 \circ
\delta_\lambda(\xi_0^{-1}\circ \xi)\in \Omega$ for $0<\lambda
<1$.}. The following are equivalent:
\begin{enumerate}
\item $u$ is $\mathcal H$--convex.
\item Given $\xi_0\in \Omega$
\begin{equation}\label{eq:convexityinh1}
u(\xi_0 \circ \delta_\lambda(\xi_0^{-1}\circ \xi)) \leq u(\xi_0)+
\lambda (u(\xi) - u(\xi_0)),
\end{equation}
for all $\xi\in \Pi_{\xi_0}$ and $0\leq \lambda \leq 1.$
\end{enumerate}
\end{proposition}

\begin{remark}\rm
From Proposition \ref{prop:H-conveximpliesconvexonlines}(2) we
have that if $u$ is convex in the standard sense, then $u$ is
$\mathcal H$--convex. However, the gauge function $\rho(x,y,t)=
\left((x^2+y^2)^2+t^2\right)^{1/4}$ is $\mathcal H$--convex but is
not convex in the standard sense, see Proposition
\ref{prop:Hofdistanceiszero}.
\end{remark}

\section{Comparison Principle}\label{sec:comparisonprinciple}

We prove the following.

\begin{theorem}\label{thm:comparisonprinciple}
Let $u,\varv\in C^2(\bar \Omega)$ such that $u+\varv$ is $\mathcal
H$--convex in $\Omega$ satisfying $\varv=u$ on $\partial \Omega$
and $\varv<u$ in $\Omega.$ Then
\[
\int_\Omega \left\{ \det \mathcal H(u)+ 12\, (\partial_t u)^2
\right\}\,dz \leq \int_\Omega \left\{ \det \mathcal H(\varv)+ 12\,
(\partial_t \varv)^2 \right\}\,dz,
\]
and
\[
\int_\Omega \trace {\mathcal H}(u) \,dz \leq \int_\Omega \trace
{\mathcal H}(\varv) \,dz.
\]
\end{theorem}

\begin{proof}
If $Z=\alpha_1\,\partial_{x_1} +\alpha_2\,\partial_{x_2}+
\alpha_3\,\partial_{x_3}$ is a smooth vector field, then
\begin{equation}\label{eq:integrationbyparts}
\int_\Omega Zu \,dx = \int_{\partial \Omega} \nu_Z \,u\,d\sigma(x)
- \int_\Omega ((\alpha_1)_{x_1}+(\alpha_2)_{x_2}
+(\alpha_3)_{x_3})\,u\,dx,
\end{equation}
where $\nu=(\nu_1, \nu_2, \nu_3)$ is the outer unit normal to
$\partial \Omega$ and $\nu_Z=\alpha_1 \nu_1 + \alpha_2 \nu_2 +
\alpha_3 \nu_3.$

Since $\varv=u$ on $\partial \Omega$, $\varv<u$ in $\Omega$ and
both functions are smooth up to the boundary, it follows that the
normal to $\partial \Omega$ is
$\nu=\dfrac{D(\varv-u)}{|D(\varv-u)|}$, and therefore $\nu_X=
\dfrac{X(\varv-u)}{|D(\varv-u)|}$ and $\nu_Y=
\dfrac{Y(\varv-u)}{|D(\varv-u)|}$.
Set
\[
S(u)
= \det \mathcal H(u)
= X^2u \, Y^2u - \left( \left( \dfrac{XY+YX}{2} \right)u\right)^2.
\]
We have
\[
\dfrac{\partial S(u)}{\partial r_{11}}=Y^2 u;
\quad
\dfrac{\partial S(u)}{\partial r_{12}}=-\left( \dfrac{XY+YX}{2} \right)u;
\]
\[
\dfrac{\partial S(u)}{\partial r_{21}}=-\left( \dfrac{XY+YX}{2} \right)u;
\quad
\dfrac{\partial S(u)}{\partial r_{22}}=X^2u.
\]
Let $0\leq s \leq 1$ and $\varphi(s)=S(s\,u+ (1-s)\,\varv)$.
Then
\begin{align*}
&\int_\Omega \{S(u)-S(\varv)\}\,dz\\
&=
\int_0^1 \int_\Omega \varphi'(s)\,dzds\\
&=
\int_0^1 \int_\Omega
\left\{  \sum_{i,j=1}^2 \frac{\partial S}{\partial r_{ij}}(\varv+s(u-\varv))\,
\left( \dfrac{X_iX_j + X_j X_i}{2}\right)(u-\varv)\right\}\,dzds\\
&=
\int_0^1 \int_\Omega
\left\{  \sum_{i,j=1}^2 \frac{\partial S}{\partial r_{ij}}(\varv+s(u-\varv))\,
(X_iX_j)(u-\varv)\right\}\,dzds\text{  since $S_{ij}$ is symmetric}\\
&=
\int_0^1 \int_\Omega
\left\{  \sum_{i,j=1}^2 X_i\left( \frac{\partial S}{\partial r_{ij}}(\varv+s(u-\varv)) \,X_j(u-\varv)\right)
-
X_i\left( \frac{\partial S}{\partial r_{ij}}(\varv+s(u-\varv))\right)
X_j(u-\varv) \right\}\,dzds\\
&=
A- B.
\end{align*}
We have
\begin{align*}
A
&=
\int_0^1\int_\Omega
X\left[ Y^2(\varv+s(u-\varv))\, X(u-\varv)\right]
-Y\left[ \left( \dfrac{XY+YX}{2} \right)(\varv+s(u-\varv))\, X(u-\varv)\right]\\
&\qquad
-X\left[ \left( \dfrac{XY+YX}{2} \right)(\varv+s(u-\varv))\, Y(u-\varv)\right]
+ Y\left[ X^2(\varv+s(u-\varv))\, Y(u-\varv)\right] \, dzds\\
&=
\int_0^1\int_{\partial\Omega} \nu_X Y^2(\varv+s(u-\varv))\, X(u-\varv) -
\nu_Y \, \left( \dfrac{XY+YX}{2}\right)(\varv+s(u-\varv))\, X(u-\varv)\\
&\qquad -
\nu_X \,\left( \dfrac{XY+YX}{2}\right)(\varv+s(u-\varv))\, Y(u-\varv) +
\nu_Y \, X^2(\varv+s(u-\varv))\, Y(u-\varv)\,dzds\\
&=
-
\int_0^1\int_{\partial\Omega}
\left\{ Y^2(\varv+s(u-\varv))\, X(u-\varv)^2 \right. \\
&\qquad \qquad +
\left( \dfrac{XY+YX}{2}\right)(\varv+s(u-\varv)) \, X(u-\varv)\, Y(u-\varv)\\
&\qquad \qquad +
\left( \dfrac{XY+YX}{2}\right)(\varv+s(u-\varv)) \, Y(u-\varv)\, X(u-\varv)\\
&\qquad \qquad -
\left. X^2(\varv+s(u-\varv))\, Y(u-\varv)^2 \right\} \dfrac{1}{|D(\varv-u)|}\,dzds\\
&=
-
\int_0^1\int_{\partial\Omega}  \left\langle \mathcal H(\varv+s(u-\varv))\,
(X(u-\varv),Y(u-\varv)),
(X(u-\varv),Y(u-\varv)) \right\rangle \dfrac{1}{|D(\varv-u)|} \,dzds\\
&\leq
0.
\end{align*}

We now calculate $B$
\begin{align*}
B
&=
\int_0^1\int_\Omega
\sum_{i,j=1}^2
X_i\left( \frac{\partial S}{\partial r_{ij}}(\varv+s(u-\varv))\right)
X_j(u-\varv)\,dzds\\
&=
\int_0^1\int_\Omega
X\left( Y^2(\varv+s(u-\varv)) \right)\,X(u-\varv)-
X\left( \dfrac{XY+YX}{2} \right)(\varv+s(u-\varv))\, Y(u-\varv)\\
&\qquad -
Y\left( \dfrac{XY+YX}{2} \right)(\varv+s(u-\varv))\, X(u-\varv)+
Y\left( X^2(\varv+s(u-\varv)) \right)\,Y(u-\varv)\,dzds\\
&=
\frac12 \int_0^1\int_\Omega
X\left( Y^2(\varv+s(u-\varv)) \right)\,X(u-\varv)-
X\left( XY+YX \right)(\varv+s(u-\varv))\, Y(u-\varv)\\
&\qquad -
Y\left( XY+YX \right)(\varv+s(u-\varv))\, X(u-\varv)+
Y\left( X^2(\varv+s(u-\varv)) \right)\,Y(u-\varv)\,dzds\\
&\qquad
+
\frac12 \int_0^1\int_\Omega
X\left( Y^2(\varv+s(u-\varv)) \right)\,X(u-\varv)
+
Y\left( X^2(\varv+s(u-\varv)) \right)\,Y(u-\varv)\,dzds\\
&=
\frac12 \int_0^1\int_\Omega
(XY-YX)Y(\varv+s(u-\varv)) \,X(u-\varv)-
(XY-YX)X(\varv+s(u-\varv))\, Y(u-\varv)\,dzds\\
&\qquad -
\frac12 \int_0^1\int_\Omega
\left( X^2Y(\varv+s(u-\varv)) \,Y(u-\varv)
+
Y^2X(\varv+s(u-\varv)) \,X(u-\varv)\right)\,dzds\\
&\qquad
+
\frac12 \int_0^1\int_\Omega
X\left( Y^2(\varv+s(u-\varv)) \right)\,X(u-\varv)
+
Y\left( X^2(\varv+s(u-\varv)) \right)\,Y(u-\varv)\,dzds\\
&=
 \int_0^1\int_\Omega
-2\partial_tY(\varv+s(u-\varv)) \,X(u-\varv)+
2\partial_t X(\varv+s(u-\varv))\, Y(u-\varv) \,dzds\\
&\qquad +
\frac12 \int_0^1\int_\Omega
\left( (XY^2 -Y^2X)(\varv+s(u-\varv)) \,X(u-\varv)
+
(YX^2-X^2Y)(\varv+s(u-\varv)) \,Y(u-\varv)\right)\,dzds
\end{align*}
On the other hand,
\[
XY^2 -Y^2X
=XY^2-YXY+YXY -Y^2X
=
(XY-YX)Y + Y(XY-YX)= -4\partial_tY -4 Y\partial_t= -8 Y\partial_t,
\]
and
\[
YX^2 -X^2Y
=YX^2-XYX+XYX -X^2Y
=
(YX-XY)X + X(YX-XY)= 4\partial_tX + 4X\partial_t= 8 X\partial_t.
\]
Therefore
\begin{align*}
B
&=
\int_0^1\int_\Omega
-2 \partial_tY(\varv+s(u-\varv)) \,X(u-\varv)+
2 \partial_t X(\varv+s(u-\varv))\, Y(u-\varv) \,dzds\\
&\qquad
\int_0^1\int_\Omega
\left( - 4Y\partial_t (\varv+s(u-\varv)) \,X(u-\varv)
+ 4
X\partial_t(\varv+s(u-\varv)) \,Y(u-\varv)\right)\,dzds\\
&=
\int_0^1\int_\Omega
\left( - 6Y\partial_t (\varv+s(u-\varv)) \,X(u-\varv)
+ 6
X\partial_t(\varv+s(u-\varv)) \,Y(u-\varv)\right)\,dzds\\
&=
6\int_0^1\int_\Omega
\partial_t (\varv + s(u-\varv))\, YX(u-\varv)\,dzds
-
6\int_0^1\int_{\partial \Omega}
\partial_t (\varv + s(u-\varv))\, X(u-v)\,\nu_Y\,d\sigma(z)\,ds\\
&-6
\int_0^1\int_\Omega
\partial_t (\varv + s(u-\varv))\, XY(u-\varv)\,dzds
+
6\int_0^1\int_{\partial \Omega}
\partial_t (\varv + s(u-\varv))\, Y(u-\varv)\,\nu_X\,d\sigma(z)\,ds\\
&=
6\int_0^1\int_\Omega
\partial_t (\varv + s(u-\varv))\, YX(u-\varv)\,dzds
-6
\int_0^1\int_\Omega
\partial_t (\varv + s(u-\varv))\, XY(u-\varv)\,dzds\\
&=
6\int_0^1\int_\Omega
\partial_t (\varv + s(u-\varv))\, (YX-XY)(u-\varv)\,dzds\\
&=
24\int_0^1\int_\Omega
\partial_t (\varv + s(u-\varv))\, \partial_t(u-\varv)\,dzds\\
&=
12\int_\Omega
\partial_t (u-\varv)\, \partial_t(u+\varv)\,dzds\\
&=
12\int_\Omega
(\partial_t u)^2 -(\partial_t\varv)^2\,dzds.
\end{align*}
This completes the proof of the theorem.

\end{proof}

\setcounter{equation}{0}
\section{Weak maximum principle}\label{sec:weakmaxprinciple}

Let $A=(a_{ij})$ be a $2\times 2$ symmetric matrix such that
$A\geq 0$, and $\trace A>0$, $a_{ij}\in C(D)$ where $D\subset
\R^3$ is an open set; $X_1=X, X_2=Y$, and $L= \sum_{i,j=1}^2
a_{ij}(\xi) X_i X_j.$
\begin{theorem}\label{thm:weakmaxprinciple}
Let $\Omega$ be a bounded open set in $\R^3$, and $w\in
C^2(\Omega)$. If $Lw\geq 0$ in $\Omega$ and $\limsup_{\xi\to
\xi_0}w(\xi)\leq 0$ for each $\xi_0\in \partial \Omega$, then
$w\leq 0$ in $\Omega.$
\end{theorem}

To prove this theorem we need two lemmas.

\begin{lemma}\label{lm:geometric}
Let $\Omega\subset \R^3$ be an open bounded set, and $w\in
C(\Omega)$. Then there exists $\xi_0\in \bar \Omega$ such that
$\sup_{\Omega\cap B(\xi_0,\rho)}w= \sup_\Omega w$ for every
$\rho>0$, where $B(\xi_0,\rho)$ is the Euclidean ball with radius
$\rho$ and center $\xi_0.$
\end{lemma}

\begin{lemma}\label{lm:barrier}
Let $\Omega$ be open and bounded. There exists a function $w_0\in
C^2(\Omega)$ such that $w_0>0$ and $Lw_0<0$ in $\Omega.$
\end{lemma}
\begin{proof}
Let $\lambda>0$ and choose $M\in \R$ such that $\sup_{\xi\in
\Omega}e^{\lambda \,x +\lambda \,y}<M;$ $\xi=(x,y,t).$ Let $w_0=
M-e^{\lambda \,x +\lambda \,y}.$ Then $w_0>0$ in $\Omega$, $X_1
w_0= - \lambda e^{\lambda \,x}$, $X_1^2 w_0= - \lambda^2
e^{\lambda \,x}$, $X_2 w_0= - \lambda e^{\lambda \,y}$, $X_2^2
w_0= - \lambda^2 e^{\lambda \,y}$, and $X_1 X_2 w_0= X_2 X_1
w_0=0.$ Hence $Lw_0=-\lambda^2 (a_{11} \,e^{\lambda \,x} +
a_{22}\, e^{\lambda \,y})<0$ in $\Omega.$
\end{proof}

\begin{proof}[Proof of Theorem \ref{thm:weakmaxprinciple}]
First assume that $Lw>0$ in $\Omega.$ By Lemma \ref{lm:geometric},
there exists $\xi_0\in \bar \Omega$ such that $\sup_{\Omega\cap
B(\xi_0,\rho)}w= \sup_\Omega w$ for every $\rho>0$. If $\xi_0\in
\Omega$, then $w(\xi_0)=\sup_\Omega w$ and so $Dw(\xi_0)=0$ and
$D^2w(\xi_0)\leq 0.$ Hence
\begin{align*}
0<Lw(\xi_0) &=
\trace \left( A \left[\begin{matrix} X^2w & XYw \\
YXw  & Y^2w
\end{matrix}\right] \right)(\xi_0)\\
&=
\trace \left( A \left[\begin{matrix} X^2w & (XYw + YXw)/2 \\
(XYw + YXw)/2  & Y^2w
\end{matrix}\right] \right)(\xi_0)\\
&= \trace \left( A \left[\begin{matrix}
1 & 0 & 2y \\
0 & 1 & -2x
\end{matrix}\right]
D^2 w \left[\begin{matrix}
1 & 0 \\
0 & 1 \\
2y & -2x
\end{matrix}\right]
\right)(\xi_0)\\
&= \trace \left( \left[\begin{matrix}
1 & 0 \\
0 & 1 \\
2y & -2x
\end{matrix}\right]
A \left[\begin{matrix}
1 & 0 & 2y \\
0 & 1 & -2x
\end{matrix}\right]
D^2 w
\right)(\xi_0)\\
&= \trace (\tilde A D^2 w)(\xi_0)\leq 0,
\end{align*}
since $\tilde A\geq 0$ and $D^2w(\xi_0)\leq 0.$ This is a
contradiction. Hence $\xi_0\in \partial \Omega$ and consequently
$w\leq 0$ in $\Omega.$ If $Lw\geq 0$ in $\Omega$, then for each
$\varepsilon>0$ we set $w_\varepsilon = w -\varepsilon \,w_0$ with
$w_0$ as in Lemma \ref{lm:barrier}. We have $Lw_\varepsilon = Lw
-\varepsilon Lw_0>0$ and $\limsup_{\xi\to \xi_0}
w_\varepsilon(\xi)\leq \limsup_{\xi\to \xi_0} w(\xi)\leq 0$ for
each $\xi_0\in \partial \Omega.$ By the previous argument,
$w_\varepsilon\leq 0$ in $\Omega$ for each $\varepsilon>0$, and so
$w\leq 0$.
\end{proof}

Let
\[
\mathcal H^*(u)=\left[\begin{matrix} Y^2u & -(XYu + YXu)/2 \\
-(XYu + YXu)/2  & X^2u
\end{matrix}\right].
\]
We have
\[
\det \mathcal H(u) = \frac12 \, \trace (\mathcal H^*(u) \,\mathcal
H(u)),
\]
and
\begin{equation}\label{eq:tracecommute}
\trace (\mathcal H^*(u) \,\mathcal H(\varv)) = \trace (\mathcal
H^*(\varv) \,\mathcal H(u)).
\end{equation}
From Theorem \ref{thm:weakmaxprinciple} we obtain the following
comparison principle.

\begin{proposition}\label{prop:comparisonprinciple}
Let $\Omega\subset \R^3$ be an open bounded set, $u,\varv\in
C^2(\Omega)$ such that $u+\varv$ is $\mathcal H$--convex, and
$\trace \{ \mathcal H(u+\varv)\}>0$. If $\det \mathcal H(u)\geq
\det \mathcal H(\varv)$ in $\Omega$ and $u\leq \varv$ on $\partial
\Omega$, then $u\leq \varv$ in $\Omega.$
\end{proposition}

\begin{proof}
We have
\begin{align*}
0&\leq \det \mathcal H(u)- \det \mathcal H(\varv)\\
&= \frac12 \, \left( \trace (\mathcal H^*(u) \,\mathcal H(u)) -
\trace (\mathcal H^*(\varv) \,\mathcal H(\varv))\right)\\ &=
\frac12 \, \left( \trace (\mathcal H^*(u) \,\mathcal H(u-\varv)) +
\trace ( (\mathcal H^*(u)- \mathcal H^*(\varv)) \,\mathcal
H(\varv))\right)\\ &= \frac12 \, \left( \trace (\mathcal H^*(u)
\,\mathcal H(u-\varv)) + \trace ( \mathcal H^*(\varv) \,\mathcal
H(u-\varv))\right) \qquad \text{by \eqref{eq:tracecommute}}\\ &=
\frac12 \,  \trace (\mathcal H^*(u+\varv) \,\mathcal H(u-\varv))\\
&= \frac12 \,  \trace (\mathcal H^*(u+\varv) \,\mathcal H(w)),
\end{align*}
where $w=u-\varv\leq 0$ on $\partial \Omega$. Applying Theorem
\ref{thm:weakmaxprinciple} to $w$ with $A=\mathcal H^*(u+\varv)$,
the proposition follows.
\end{proof}

\subsection{A comparison Principle}
As a consequence of Proposition \ref{prop:comparisonprinciple} we
get that "cones" that agreeing with an $\mathcal H$--convex
function $u$ on the boundary of a ball $B$ are above $u$ inside
$B$.

\begin{proposition}\label{prop:Hofdistanceiszero}
Let $d(\xi,\xi_0)=\|\xi_0^{-1}\circ \xi \|$, $\|\xi\|= ((x^2 +
y^2)^2 + t^2)^{1/4}$, $\xi=(x,y,t)$, $\Omega= \{\xi\in \R^3:
0<d(\xi,\xi_0) <R\}$, and $\varv(\xi)= m
\left(\dfrac{d(\xi,\xi_0)}{R} -1\right)$. If $m\geq 0$, then
$\varv$ is $\mathcal H$--convex in $\Omega$, $\det \mathcal
H(\varv)=0$ in $\Omega$, and $\det \mathcal H(\varv)$ is
integrable in $\bar \Omega$.
\end{proposition}

\begin{proof}
If $\zeta\in \R^3$ and $g(\xi)=f(\zeta\circ \xi)$, then
$Xg(\xi)=(Xf)(\zeta\circ \xi)$ and $Yg(\xi)=(Yf)(\zeta\circ \xi)$.
Therefore we can assume that $\xi_0=0.$ Let $r=(x^2+y^2)^2+t^2$
and $h\in C^1((0,+\infty))$. Then $Xr= 4x^3 + 4xy^2+4yt$,
$Yr=4yx^2+ 4 y^3 -4xt$, $X^2r=Y^2r=12(x^2+y^2)$, $YXr=4t$, and
$XYr=-4t$. If $u(x,y,t)=h(r)$, then $Xu=h'(r)\,Xr$,
$Yu=h'(r)\,Yr$, $X^2u=h''(r)\,(Xr)^2+h'(r)\, X^2r$,
$Y^2u=h''(r)\,(Yr)^2+h'(r)\, Y^2r$, $XYu=h''(r)\,Xr Yr +h'(r)\,
XYr$, $YXu=h''(r)\,Yr Xr +h'(r)\, YXr$. Thus
\begin{equation}\label{eq:formulafordetH(u)}
\det \mathcal H(u)=48\,(x^2+y^2)^2 \,\{4\,r \,h''(r)+ 3
\,h'(r)\}\, h'(r).
\end{equation}
Therefore $\det \mathcal H(u)=0$ if $h'(r)=0$ or $4\,r \,h''(r)+ 3
\,h'(r)=0$, that is, $h(r)=C$ or $h(r)=r^{1/4}$. If
$h(r)=r^{1/4}$, then $X^2h(r)= 3\, r^{-7/4}\,
(y(x^2+y^2)-xt)^2\geq 0$ and $Y^2h(r)= 3\, r^{-7/4}\,
(x(x^2+y^2)+yt)^2\geq 0$, and o $r^{1/4}$ is $\mathcal H$--convex
in $\R^3\setminus \{0\}$.

On the other hand, $\det \mathcal H(u)\leq C\, r^{-1/2}$ and so
$\int_{r^{1/4}\leq R} \det \mathcal H(u)\,dz \leq C\,
\int_{r^{1/4}\leq R} r^{-1/2}\,dz= C\,\int_0^R
\rho^{Q-1}\,\rho^{-2}\,d\rho=C\,R^2$, since $Q=4$.
\end{proof}

\begin{proposition}\label{prop:ulessthandistance}
Let $u\in C^2(\Omega)$ be $\mathcal H$--convex, with $\Omega=
\{\xi\in \R^3: 0<d(\xi,\xi_0) <R\}$, and $u\leq 0$ on $\{\xi\in
\R^3: d(\xi,\xi_0)=R\}$. Then $u\leq \varv$, where $\varv$ is
defined in Proposition \ref{prop:Hofdistanceiszero} with
$m=-u(\xi_0).$
\end{proposition}

\begin{proof}
Let $\varepsilon>0$, $\xi_0=(x_0,y_0,t_0)$, $\xi=(x,y,t)$,
\[
u_\varepsilon (\xi)=u(\xi)+ \varepsilon \, (x^2+y^2),
\]
and
\[
\varv_\varepsilon (\xi) = -(1-\sqrt{\varepsilon})\, u(\xi_0)\,
\left(\dfrac{d(\xi,\xi_0)}{(1-\sqrt{\varepsilon})\,R} -1 \right).
\]
We first claim that $u_\varepsilon(\xi)\leq
\varv_\varepsilon(\xi)$ for all $\xi\in \partial \Omega$ and for
all $\varepsilon$ sufficiently small. Indeed, if $\xi=\xi_0$, then
$u_\varepsilon(\xi_0)\leq \varv_\varepsilon(\xi_0)$ if and only if
$\sqrt{\varepsilon} \, (x_0^2 + y_0^2)\leq -u(\xi_0)$ which holds
for all $\varepsilon$ sufficiently small. On the other hand, if
$d(\xi,\xi_0)=R$, then
$\varv_\varepsilon(\xi)=-\sqrt{\varepsilon}\,u(\xi_0)$ and
$u_\varepsilon(\xi) \leq \varepsilon \, (x^2+y^2)\leq \varepsilon
\,\max_{d(\xi,\xi_0)=R} (x^2+y^2) =\varepsilon \,M.$ Hence
$u_\varepsilon(\xi)\leq \varv_\varepsilon(\xi)$ on
$d(\xi,\xi_0)=R$ if $\sqrt{\varepsilon}\, M \leq -u(\xi_0)$ which
again holds for all $\varepsilon$ sufficiently small.

We also have
\begin{equation}\label{eq:detofuepsilon}
\det \mathcal H(u_\varepsilon)= \det \mathcal H(u) + 2 \varepsilon
\,\trace \mathcal H(u) + 4\,\varepsilon^2> 0=\det \mathcal
H(\varv_\varepsilon)
\end{equation}
in $\Omega,$ and $ \trace \{\mathcal H(u_\varepsilon +
\varv_\varepsilon)  \}= \trace \mathcal H(u)+ 8\,\varepsilon +
\trace \mathcal H(\varv_\varepsilon)>0. $ Therefore from
Proposition \ref{prop:comparisonprinciple} we get
$u_\varepsilon\leq \varv_\varepsilon$ in $\Omega,$ and the
proposition follows letting $\varepsilon\to 0.$

\end{proof}

As a consequence of these propositions we get the following
extension of Theorem \ref{thm:comparisonprinciple} needed in the
proof of the maximum principle Theorem \ref{ABP}.
\begin{theorem}\label{thm:comparisonprincipleperforated}
Let $\Omega= \{\xi\in \R^3: 0<d(\xi,\xi_0) <R\}$, and let
$\varv\in C^2(\bar B_R(\xi_0))$ be $\mathcal H$--convex in
$\Omega$ satisfying $\varv=0$ on $\partial B_R(\xi_0)$ and set
$u(\xi) = -\varv(\xi_0)\, \left(\dfrac{d(\xi,\xi_0)}{R}
-1\right)$. Then
\[
\int_{B_R(\xi_0)} \left\{ \det \mathcal H(u)+ 12\, (\partial_t
u)^2 \right\}\,d\xi \leq \int_{B_R(\xi_0)} \left\{ \det \mathcal
H(\varv)+ 12\, (\partial_t \varv)^2 \right\}\,d\xi.
\]
\end{theorem}

\begin{proof}
From Proposition \ref{prop:ulessthandistance} we have that
$\varv\leq u$ in $B_R(\xi_0)$. Let $\varepsilon >0$, we claim that
\begin{align}\label{eq:boundplusOepsilon}
&\int_{B_R(\xi_0)\setminus B_\varepsilon(\xi_0)} \left\{ \det
\mathcal H(u)+ 12\, (\partial_t u)^2
\right\}\,d\xi\notag \\
&\qquad \leq \int_{B_R(\xi_0)\setminus B_\varepsilon(\xi_0)}
\left\{ \det \mathcal H(\varv)+ 12\, (\partial_t \varv)^2
\right\}\,d\xi + O(\varepsilon^{1/4}),
\end{align}
as $\varepsilon\to 0.$ We may assume by the invariance of the
vector fields that $\xi_0=0$. Since the functions $u,\varv$ are
both convex and $C^2$ except at $0$, we proceed as in the proof of
Theorem \ref{thm:comparisonprinciple} applied to the open set
$\Omega_\varepsilon=B_R(0)\setminus B_\varepsilon(0)$. The sum of
the integrals $I,II,III$ and $IV$ contains now the boundary terms
\begin{align*}
&-\int_{d(\xi)=\varepsilon} Y^2 (u + \varv) \, X(u-\varv) \,
\dfrac{Xd}{|Dd|}\,d\sigma(\xi) - \int_{d(\xi)=\varepsilon} X^2 (u
+ \varv) \, Y(u-\varv) \,
\dfrac{Yd}{|Dd|}\,d\sigma(\xi)\\
&- \int_{d(\xi)=\varepsilon} XY (u + \varv) \, X(u-\varv) \,
\dfrac{Yd}{|Dd|}\,d\sigma(\xi) - \int_{d(\xi)=\varepsilon} YX (u +
\varv) \, Y(u-\varv) \, \dfrac{Xd}{|Dd|}\,d\sigma(\xi),
\end{align*}
where $d(\xi)=d(\xi,0).$ We shall prove that each summand is
$O(\varepsilon^{1/4}).$ Each of these summands basically have the
same behaviour as $\varepsilon\to 0.$ Using the computations used
in the proof of Proposition \ref{prop:Hofdistanceiszero}, we see
for example that
\begin{align}
&J=\int_{d(\xi)=\varepsilon} Y^2 (u + \varv) \, X(u-\varv) \,
\dfrac{Xd}{|Dd|}\,d\sigma(\xi)\notag\\
&\qquad \leq C\, \int_{d(\xi)=\varepsilon}
r^{-11/4}\,(x(x^2+y^2)+yt)^2\,
(4x^3 + 4 xy^2 + 4 yt)^2\,\dfrac{d\sigma(\xi)}{|Dd|}\notag\\
&\qquad \leq C\,\varepsilon^{-11/4}\,\int_{d(\xi)=\varepsilon}
\dfrac{d\sigma(\xi)}{|Dd|}\label{eq:valueofintegralcoarea}.
\end{align}
On the other hand, from the coarea formula
\[
\int_0^t \int_{d(\xi)=s} \dfrac{d\sigma(\xi)}{|Dd|} \,ds=
\int_{d(\xi)\leq t}d\xi= C\, t^4.
\]
So $\displaystyle \int_{d(\xi)=s}
\dfrac{d\sigma(\xi)}{|Dd|}=C\,s^3$ and inserting this value in
\eqref{eq:valueofintegralcoarea} we obtain that
$J=O(\varepsilon^{1/4}).$

Following the method of proof of Theorem
\ref{thm:comparisonprinciple} we integrate by parts once again and
we now obtain the boundary terms
\[
\int_{d(\xi)=\varepsilon} X(u + v) \, \partial_t(u-v) \,
\dfrac{Yd}{|Dd|}\,d\sigma(\xi),
\]
and
\[
\int_{d(\xi)=\varepsilon} Y(u + v) \, \partial_t(u-v) \,
\dfrac{Xd}{|Dd|}\,d\sigma(\xi).
\]
These integrals can be handled as before obtaining again that they
are $O(\varepsilon^{1/4}).$ Therefore \eqref{eq:boundplusOepsilon}
holds and the theorem follows letting $\varepsilon\to 0.$
\end{proof}

As a consequence of Proposition \ref{prop:comparisonprinciple} we
obtain that $\mathcal H$--convex functions are Lipschitz with
respect to the distance $d$.

\begin{proposition}\label{prop:uHconveximplieslipschitz}
Let $\Omega\subset \R^3$ be an open set and $u\in C(\Omega)$
$\mathcal H$--convex in $\Omega$. Then for each ball $\bar
B\subset \Omega$ there exists a constant $C_B$ such that
$|u(x)-u(y)|\leq C_B\, d(x,y)$ for all $x,y\in B$.
\end{proposition}

\begin{proof}
We can assume that $u\in C^2(\Omega)$ and let $B_d(x_0,2R)\subset
\Omega$. Let $y\in B_d(x_0,R)$ and $\phi(x)=u(x)-u(y)+\varepsilon
\left( (x_1-y_1)^2+ (x_2-y_2)^2\right)$; $x=(x_1,x_2,x_3)$,
$y=(y_1,y_2,y_3)$, with $x\in B_d(y,R)$. We have $\mathcal H(\phi
+ C_\e \,d(\cdot ,y))>0$ in $B_d(y,R)$ and $\phi(x)\leq C_\e \,
d(x,y)$ for $d(x,y)=R$ where $C_\e=\dfrac{\osc_{B_d(x_0,2R)} u \,+
\varepsilon \, \diam(B_d(x_0,2R))^2}{R}$. We have $\det \mathcal
H(\phi)\geq \det \mathcal H(d(\cdot,y))$ in $B_d(y,R)\setminus
\{y\}$ so by the comparison principle Proposition
\ref{prop:comparisonprinciple} we get that $\phi(x)\leq C_\e\,
d(x,y)$ for $x\in B_d(y,R)$. Letting $\e\to 0$ we get
$u(x)-u(y)\leq C\, d(x,y)$ for $x\in B_d(y,R)$ with
$C=\dfrac{\osc_{B_d(x_0,2R)} u}{R}$ and $y\in B_d(x_0,R)$. If
$x,y\in B_d(x_0,R/4)$, then $x\in B_d(y,R/2)$ and so $y\in
B_d(x,R)$ and by the previous inequality we get $u(y)-u(x)\leq C\,
d(y,x)=C\, d(x,y)$. Therefore we obtain $|u(x)-u(y)|\leq C\,
d(x,y)$ for all $x,y\in B_d(x_0,R/4)$.
\end{proof}

\section{Maximum Principle}\label{sec:maximumprinciple}


\begin{proposition}\label{core}
Let $u$ be ${\mathcal H}$--convex in $\Omega$ open and bounded.
Suppose $u\leq 0$ on $\partial \Omega$. Then $u\leq 0$ in
$\Omega$.
\end{proposition}
\begin{proof}
Let $\varepsilon>0$ and $u_\varepsilon(x,y,t)=
u(x,y,t)+\varepsilon \,(x^2 +y^2).$ We have $\mathcal
H(u_\varepsilon)= \mathcal H(u) + 2\varepsilon \,\text{Id}$, so
\[
\det \mathcal H(u_\varepsilon)= \det \mathcal H(u) + 2 \varepsilon
\,\trace \mathcal H(u) + 4\,\varepsilon^2.
\]
Since $\det \mathcal H(\varepsilon(x^2 + y^2))=4\,\varepsilon^2$,
we get $\det \mathcal H(u_\varepsilon)\geq \det \mathcal
H(\varepsilon(x^2 + y^2))$ in $\Omega$. Also $u_\varepsilon \leq
\varepsilon(x^2 + y^2)$ on $\partial \Omega$, and $\trace
\{\mathcal H(u_\varepsilon + \varepsilon (x^2 + y^2)  \}= \trace
\mathcal H(u)+ 8\,\varepsilon >0$. The proposition then follows
from Proposition \ref{prop:comparisonprinciple}.
\end{proof}

\begin{proposition}\label{core1}
Let $u$ be ${\mathcal H}$--convex in $\Omega$ open and bounded.
Suppose $u\leq 0$ on $\partial \Omega$. Then $u\leq 0$ in
$\Omega$. Moreover, if there is $\xi_0\in \Omega$ such that
$u(\xi_0)=0$ then $u\equiv 0$ in $\Omega.$
\end{proposition}
\begin{proof}
Define
\[
L := X^2+Y^2
\]
the Kohn Laplacian on the Heisenberg group. Since $u$ is
${\mathcal H}$--convex in $\Omega$ then $\trace \mathcal Hu
=Lu\geq 0.$ Hence, by the maximum principle for $L,$ we get $u\leq
0$ in $\Omega.$ Moreover, if there is $\xi_0\in \Omega$ such that
$u(\xi_0)=0$ then $u$ has a maximum at an interior point and by
strong maximum principle for $L$, see \cite{Bony:maxpciple}, we
get $u\equiv 0$ in $\Omega.$
\end{proof}

The following lemma will be used repeatedly in the proof of
Proposition \ref{esti}.
\begin{lemma}\label{lm:mainlemma}
Let $\xi_0 \in B_R(0)$ and $\xi\in \Pi_{\xi_0} \cap B_R(0)$. Let
$\lambda>0$ be such that
\[
\xi'=\xi_0 \circ \delta_\lambda (\xi_0^{-1}\circ \xi) \in
\Pi_{\xi_0} \cap \partial B_R(0).
\]
Suppose $u$ is $\mathcal H$--convex in $B_R(0)$ and $u=0$ on
$\partial B_R(0)$. Then:
\begin{enumerate}
\item If
$\xi_0=(x_0,y_0,t_0)$ and $\xi=(0,0,t_0)$, then $\lambda\geq 2$
and
\begin{equation}\label{eq:pointwiseestimateofu}
u(\xi)\leq \dfrac{1}{2}\, u(\xi_0).
\end{equation}

\item If $0<\alpha,\beta <1$, $\alpha +\beta <1$,
$\rho(\xi_0)\leq \alpha \, R$ and $d(\xi_0,\xi)\leq \beta \, R$,
then $\lambda \geq \dfrac{1-\alpha}{\beta}$ and
\begin{equation}\label{eq:pointwiseestimateofubis}
u(\xi)\leq \dfrac{1-\alpha -\beta }{1-\alpha}\, u(\xi_0).
\end{equation}
\end{enumerate}
\end{lemma}

\begin{proof}
To prove the first part of (1), if $\eta=(x,y,t)\in \Pi_{\xi_0}$,
then we have that
\[
\xi_0\circ \delta_\lambda (\xi_0^{-1}\circ \eta)= (x_0+ \lambda
(x-x_0), y_0+\lambda (y-y_0), t_0+\lambda (t-t_0)),
\]
in particular, $\xi'=((1-\lambda)x_0, (1-\lambda)y_0, t_0).$ Hence
\begin{align*}
R^4 &= \rho(\xi')^4 = \left( (1-\lambda)^2 x_0^2 + (1-\lambda)^2
y_0^2 \right)^2 +
t_0^2\\
&= (1-\lambda)^4\,\left( x_0^2 +  y_0^2 \right)^2
+ \rho(\xi_0)^4 - (x_0^2 + y_0^2)^2\\
&\leq \left( (1-\lambda)^4 -1\right)\,\left( x_0^2 +  y_0^2
\right)^2 + R^4,
\end{align*}
and so $|1-\lambda|\geq 1$. Since $\lambda>0$, it follows that
$\lambda \geq 2.$

To prove the first part of (2) we write
\begin{align*}
R&=\rho(\xi') = \rho((\xi_0^{-1})^{-1}\circ \delta_\lambda
(\xi_0^{-1}\circ \xi)) \leq
\rho(\xi_0^{-1})+ \rho(\delta_\lambda (\xi_0^{-1}\circ \xi))\\
&= \rho(\xi_0)+ \lambda \, \rho(\xi_0^{-1}\circ \xi) =
\rho(\xi_0) + \lambda \, d(\xi_0,\xi)\\
& \leq \alpha \,R + \lambda \, \beta \, R,
\end{align*}
and so $\lambda \geq \dfrac{1-\alpha}{\beta}$ .

To prove \eqref{eq:pointwiseestimateofu} and
\eqref{eq:pointwiseestimateofubis}, by definition of $\xi'$ we
have that $\xi=\xi_0\circ \delta_{1/\lambda}(\xi_0^{-1}\circ
\xi')$. From \eqref{eq:convexityinh1} and since $u(\xi')=0$, it
follows that $u(\xi)\leq \left(1- \dfrac{1}{\lambda}\right)
u(\xi_0)$. Thus \eqref{eq:pointwiseestimateofu} and
\eqref{eq:pointwiseestimateofubis} follow since $u\leq 0$ in
$B_R(0).$
\end{proof}

\begin{proposition}\label{esti}
Let $u$ be ${\mathcal H}$--convex and $u=0$ on $\partial B_R(0).$
Given $\xi_0\in B_{R}(0)$ there exists a positive constant $c<1,$
depending on $d(\xi_0,\p B_R(0)),$ such that
\[
u(0)\leq c\, u(\xi_0).
\]
\end{proposition}

\begin{proof}
Let $\xi_0=(x_0,y_0,t_0)$ and
$\xi_1=\exp(-x_0X-y_0Y)(\xi_0)=(0,0,t_0)\in \Pi_{\xi_0}.$ We
obviously have that $ d(\xi_1,\xi_0)=\sqrt{x_0^2+y_0^2}\leq
d(0,\xi_0)<R. $ Applying Lemma \ref{lm:mainlemma}(1) with
$\xi_0\rightsquigarrow \xi_0$ and $\xi\rightsquigarrow \xi_1$ we
get that
\begin{equation}\label{eqes}
u(\xi_1)\leq \frac{1}{2} u(\xi_0).
\end{equation}

We shall prove that there exists a constant $C_1>0$ depending only
of the distance from $\xi_1$ to $\partial B_R(0)$ such that
\begin{equation}\label{eq:estimateu(0)}
u(0)\leq C_1\, u(\xi_1).
\end{equation}
To prove \eqref{eq:estimateu(0)} we may assume $\xi_1\ne 0,$ and
consider two cases.

{\bf Case 1.} $d(\xi_1,0)=|t_0|^{1/2}\leq R/2$.

If $t_0>0,$ define $\sigma=\dfrac{\sqrt{t_0}}{2}$ and put
\[
\begin{split}
\xi_2=&\exp(\sigma X)\xi_1=(\s,0,t_0),\\
\xi_3=&\exp(\sigma Y)\xi_2=(\s,\s,t_0-2\s^2),\\
\xi_4=&\exp(-\sigma X)\xi_3=(0,\s,t_0-2\s^2-2\s^2)
=(0,\s,t_0-4\s^2)=(0,\s,0).
\end{split}
\]
By our choice of $\s$ we have
\[
\exp(-\sigma Y)\xi_4=(0,0,t_0-4\s^2)=0.
\]
Let us remark that
\[
\s= \frac{1}{2}d(\xi_1,0)\leq R/4.
\]
We have
\[
d(\xi_1,\xi_2)= d(\xi_2,\xi_3)= d(\xi_3,\xi_4)=\sigma;
\]
\[
\rho(\xi_2)=17^{1/4}\,\sigma; \qquad \rho(\xi_3)=8^{1/4}\,\sigma;
\qquad \rho(\xi_4)=\sigma.
\]
Hence $\xi_2,\xi_3,\xi_4 \in B_R.$ Applying Lemma
\ref{lm:mainlemma}(2) with $\xi_0\rightsquigarrow \xi_1$,
$\xi\rightsquigarrow \xi_2$, $\alpha=1/2$, and $\beta=1/4$ we get
that
\[
u(\xi_2)\leq \frac{1}{2}u(\xi_1).
\]
Next, applying Lemma \ref{lm:mainlemma}(2) with
$\xi_0\rightsquigarrow \xi_2$, $\xi\rightsquigarrow \xi_3$,
$\alpha=17^{1/4}/4$, and $\beta=1/4$, we get that
\[
u(\xi_3)\leq \dfrac{3-17^{1/4}}{4-17^{1/4}} \, u(\xi_2)
<\dfrac{3}{8} \,\,u(\xi_2).
\]
Applying once again Lemma \ref{lm:mainlemma}(2) now with
$\xi_0\rightsquigarrow \xi_3$ and $\xi\rightsquigarrow \xi_4$,
$\alpha=8^{1/4}/4$, $\beta=1/4$, we get that
\[
u(\xi_4)\leq  \dfrac{3-8^{1/4}}{4-8^{1/4}}\, u(\xi_3) <
\dfrac12\,\,u(\xi_3).
\]
Define
\[
\xi^{(4)}=\xi_4\circ \delta_\l(\xi_4^{-1})\in \Pi_{\xi_4}
\]
and choose $\l>0$ such that $\xi^{(4)}\in \p B_R.$ Applying Lemma
\ref{lm:mainlemma}(2) now with $\xi_0\rightsquigarrow \xi_4$ and
$\xi\rightsquigarrow 0$, $\alpha=1/4$, $\beta=1/4$, we get that
\[
u(0)\leq \dfrac{2}{3} \, u(\xi_4).
\]
This completes the proof of \eqref{eq:estimateu(0)} for $t_0>0$.

If $t_0<0,$ define $\sigma=\dfrac{\sqrt{-t_0}}{2}$ and put
\[
\begin{split}
\xi_2=&\exp(\sigma Y)\xi_1=(0,\s,t_0),\\
\xi_3=&\exp(\sigma X)\xi_2=(\s,\s,t_0+2\s^2 ),\\
\xi_4=&\exp(-\sigma Y)\xi_3=(\s,0,t_0 +4\s^2) .
\end{split}
\]
By our choice of $\s$ we have
\[
\exp(-\sigma X)\xi_4=(0,0,t_0+4\s^2)=0.
\]
Then, arguing as in case $t_0>0$, we get \eqref{eq:estimateu(0)}.

{\bf Case 2.} $R/2<d(\xi_1,0)=|t_0|^{1/2}<R.$

Define
\begin{equation}\label{def:defofd}
d:=\frac{d(\xi_1,\p
B_R)}{\sqrt{6}}=\frac{\sqrt{R^2-|t_0|}}{\sqrt{6}}.
\end{equation}
Obviously $d^2< R^2/8.$ It is not restrictive to assume $t_0>0.$
We first prove that there exists a universal constant $0<C_2<1$
such that
\begin{equation}\label{d}
u(0,0,t_0-4d^2)\leq C_2\,u(\xi_1).
\end{equation}
Let
\begin{align*}
\xi_1&=(0,0,t_0)\\
\xi_2&=\exp(dX)(\xi_1)=(d,0,t_0)\\
\xi_3&=\exp(dY)(\xi_2)=(d,d,t_0-2d^2)\\
\xi_4&=\exp(-dX)(\xi_3)=(0,d,t_0-4d^2)\\
\xi_5&=\exp(-d Y)(\xi_4)=(0,0,t_0-4d^2).
\end{align*}
We have $\xi_{i+1}\in \Pi_{\xi_i}$ for $i=1,2,3,4.$ Let
\[
\xi_2^{(1)}=\exp(\l d X)(\xi_1)=(\l d,0,t_0) = \xi_1\circ
\delta_\lambda(\xi_1^{-1}\circ \xi_2),
\]
with $\l>0$ such that $\xi_2^1\in \Pi_{\xi_1}\cap \p B_R.$ Then
\begin{equation*}
R^4=\rho(\xi_2^{(1)})=\l^4d^4+t_0^2=\l^4d^4+(R^2-6d^2)^2=(\l^4+36)d^4+R^4-12d^2R^2,
\end{equation*}
and so
\[
12R^2=(\l^4+36)d^2\leq (\l^4+36)R^2/8
\]
which yields $\l >2.$ Hence,
\[
u(\xi_2)\leq (1/2) u(\xi_1).
\]

We have
\[
\begin{split}
\rho(\xi_2)^4&=d^4+t_0^2=d^4+(R^2-6d^2)^2\\
&=37d^4+R^4-12R^2d^2
=d^2(37d^2-12R^2)+R^4\\
&\leq d^2(37/8-12)R^2+R^4 = \left( \dfrac{1}{8} \left(
\dfrac{37}{8}-12 \right)+1 \right)\,R^4<R^4,
\end{split}
\]
and
\[
d(\xi_2,\xi_3)=d\leq \dfrac{1}{\sqrt{8}}\, R.
\]
If
\[
\xi_3^{(2)}=\exp(\l d Y)(\xi_2)=(d,\l d,t_0-2\l d^2) = \xi_2\circ
\delta_\lambda (\xi_2^{-2}\circ \xi_3)
\]
and we pick $\l>0$ such that $\xi_3^2\in \Pi_{\xi_2}\cap \p B_R,$
then applying Lemma \ref{lm:mainlemma}(2) with
$\xi_0\rightsquigarrow \xi_2$, $\xi\rightsquigarrow \xi_3$,
$\alpha=\sqrt[4]{ \dfrac{1}{8} \left( \dfrac{37}{8}-12
\right)+1}$, and $\beta=\dfrac{1}{\sqrt{8}}$, we get that
\[
u(\xi_3)\leq \dfrac{\sqrt{8}-\sqrt[4]{5}-1}{\sqrt{8}-\sqrt[4]{5}}
\, u(\xi_2) < \dfrac15 \,\, u(\xi_2).
\]

Next,
\[
\begin{split}
\rho(\xi_3)^4&=(2d^2)^2+(t_0-2d^2)^2=4d^4+(R^2-8d^2)^2=68d^4-16R^2d^2+R^4\\
&\leq d^2R^2(68/8-16)+R^4 \leq \left( \dfrac{1}{8} \left(
\dfrac{68}{8}-16 \right)+1 \right) \, R^4<R^4,
\end{split}
\]
and
\[
d(\xi_3,\xi_4)= d\leq \dfrac{1}{\sqrt{8}}\, R.
\]
Let
\[
\xi_4^{(3)}=\exp(-\l d X)(\xi_3)=((1-\l)d,d,t_0-2d^2-2\l d^2) =
\xi_3\circ \delta_\lambda(\xi_3^{-1}\circ \xi_4),
\]
with $\l>0$ such that $\xi_4^{(3)}\in \p B_R\cap \Pi_{\xi_3}.$
Applying Lemma \ref{lm:mainlemma}(2) with $\xi_0\rightsquigarrow
\xi_3$, $\xi\rightsquigarrow \xi_4$, $\alpha=\sqrt[4]{\dfrac{1}{8}
\left( \dfrac{68}{8}-16 \right)+1}=\dfrac{1}{2}$, and
$\beta=\dfrac{1}{\sqrt{8}}$, we get that
\[
u(\xi_4)< \dfrac14 u(\xi_3).
\]

We have
\[
\begin{split}
\rho(\xi_4)^4&=d^4+(t_0-4d^2)^2=d^4+(R^2-10d^2)^2=101d^4+R^4-20R^2d^2\\
&=(101d^2-20R^2)d^2+R^4\leq (101/8-20)R^2d^2+R^4\\
&\leq \left( \dfrac{1}{8} \left( \dfrac{101}{8}-20 \right)+1
\right)\,R^4<R^4,
\end{split}
\]
and
\[
d(\xi_4,\xi_5)=d\leq \dfrac{1}{\sqrt{8}}\, R.
\]
Letting
\[
\xi_5^{(4)}=\exp(-\l d Y)(\xi_4)=(0,(1-\l)d,t_0-4d^2) = \xi_4\circ
\delta_\lambda (\xi_4^{-1}\circ \xi_5)
\]
with $\l>1$ such that $\xi_5^{(4)} \in \Pi_{\xi_4}\cap \p B_R,$
and applying Lemma \ref{lm:mainlemma}(2) with
$\xi_0\rightsquigarrow \xi_4$, $\xi\rightsquigarrow \xi_5$,
$\alpha=\sqrt[4]{\dfrac{1}{8} \left( \dfrac{101}{8}-20 \right)+1}=
\dfrac{\sqrt[4]{5}}{\sqrt{8}}$, and $\beta=\dfrac{1}{\sqrt{8}}$,
we get that
\[
u(\xi_5)\leq \dfrac15 u(\xi_4).
\]

Thus, inequality \eqref{d} follows.

We now iterate the inequality \eqref{d}. Let $d_0=d$ (defined in
\ref{def:defofd}), $t_1=t_0 - 4\, d_0^2$, and in general
\[
t_{j+1}= t_j - 4\,d_j^2, \qquad \text{and} \qquad
d_j^2=\dfrac{R^2-t_j}{6}.
\]
We have $d_{j+1}^2=\dfrac{R^2-t_{j+1}}{6}=\dfrac{R^2-t_j+4\,
d_j^2}{6} = \left( 1+\dfrac23 \right)\,d_j^2$. Thus,
\begin{align}\label{eq:definitionoftj}
t_{N+1}=t_0 - 4\, \sum_{j=0}^N d_j^2 &=
t_0 - 4\,d_0^2 \sum_{j=0}^N \left( 1+\dfrac23 \right)^j\notag\\
&= t_0 - (R^2-t_0)\left( \left( 1+\dfrac23 \right)^{N+1} -1
\right).
\end{align}
Pick $N$ such that
\[
t_N\leq \dfrac{R^2}{4} < t_{N-1},
\]
which amounts
\begin{equation}\label{eq:chooseN}
N-1< \ln \left[\dfrac{3 \, R^2}{4(R^2-t_0)}\right]^{1/\ln(1+2/3)}
\leq N.
\end{equation}
We have $t_N< t_{N-1}<\cdots <t_1<t_0$ and it is easy to check
from \eqref{eq:definitionoftj}, the choice of $N$ and
\ref{def:defofd} that $t_N\geq -R^2/4.$ Therefore $(0,0,t_j)\in
B_R(0)\setminus B_{R/2}(0)$ for $0\leq j\leq N-1$ and
$(0,0,t_N)\in B_{R/2}(0)$. Iterating \eqref{d} $N$ times, then
yields
\[
u(0,0,t_N)\leq C_1^N \, u(\xi_1).
\]
Since $0<C_1<1$, there is $\gamma>0$ such that $C_1=e^{-\gamma}$,
and from \eqref{eq:chooseN} we obtain
\begin{align*}
u(0,0,t_N)&\leq C_1\,\exp\left(-\gamma \ln \left[
\dfrac{3 \, R^2}{4(R^2-t_0)}\right]^{1/\ln(1+2/3)}\right) \, u(\xi_1)\\
&=C_1\, \left[ \dfrac{4(R^2-t_0)}{3 \,
R^2}\right]^{\gamma/\ln(1+2/3)}\, u(\xi_1).
\end{align*}
Since $(0,0,t_N)\in B_{R/2}(0)$, we can apply
\eqref{eq:estimateu(0)} to get $u(0)\leq C_1 \, u(0,0,t_N)$.
Consequently,
\[
u(0)\leq C_1^2\, \left[ \dfrac{4(R^2-t_0)}{3 \,
R^2}\right]^{\gamma/\ln(1+2/3)}\, u(\xi_1),
\]
which completes the proof of \eqref{eq:estimateu(0)} in Case 2.

Finally, combining \eqref{eqes} and \eqref{eq:estimateu(0)} we
obtain the proposition.
\end{proof}

\begin{theorem}\label{ABP}
Let $u\in C^2(B_R)$ be ${\mathcal H}$--convex, $u=0$ on $\p B_R$.
If
\[
u(\xi_0)=\min_{B_R}u,
\]
then there exists a positive constant $c$, depending on
$d(\xi_0,\p B_R)$, such that
\[
|u(\xi_0)|^2\leq c \int_{B_R}(\det \mathcal H(u)+12\,u_t^2)\,dz.
\]
\end{theorem}

\begin{proof}
Define
\[
u(0)=-m
\]
and
\[
\varv(\zeta)=m\,\left(\frac{d(\zeta,0)}{R}-1\right).
\]
We have $\varv=u=0$ on $\p {B}_R$, $\varv$ is ${\mathcal
H}$--convex in ${B}_R$ and $\varv\geq u$ in ${B}_R.$ From the
comparison principle, Theorem
\ref{thm:comparisonprincipleperforated}, we then get
\[
\int_{B_R}\{\det \mathcal H(\varv)+12\,\varv_t^2 \}\,dz\leq
\int_{B_R} \{\det \mathcal H(u)+12\,u_t^2\}\,dz.
\]
Moreover,
\[
\begin{split}
\int_{B_R} \{\det \mathcal H(\varv)+12\,\varv_t^2\}\,dz
&=\left(\frac{m}{R}\right)^2
\int_{B_R}\{\det \mathcal H(d(\zeta,0))+12\,(\p_td(\zeta,0))^2\}\,d\zeta\\
&=12\left(\frac{m}{R}\right)^2R^2
\int_{B_1}(\p_td(\zeta,0))^2d\zeta\\ &= c_1 m^2
\end{split}
\]
with
\[
c_1=12\int_{B_1}(\p_td(\zeta,0))^2d\zeta>0.
\]

Let
\[
u(\xi_0)=\min_{B_R}u=-m_0.
\]
By Proposition \ref{esti} there exists a constant $0<c_2<1$ such
that
\[
m_0 \leq \frac{1}{c_2} m .
\]
Hence,
\[
m_0^2\leq  \frac{1}{c_2^2} \, m^2\leq
\frac{c_1}{c_2^2}\int_{B_R}\{\det \mathcal H(u)+12\,u_t^2\}\,dz.
\]
\end{proof}

\section{${\mathcal H}$--Measures}\label{sec:Hmeasures}

\subsection{Oscillation estimate}
In this section we prove that if $u$ is ${\mathcal H}$--convex, we
can control the integral of $\det \mathcal H(u)+12(u_t)^2$ locally
in terms of the oscillation of $u.$

Let us start with a lemma on ${\mathcal H}$--convex functions,
which is similar to the Euclidean one for convex functions.
\begin{lemma}\label{lem:comp}
If $u_1,u_2\in C^2(\Omega)$ are ${\mathcal H}$--convex, and $f$ is
convex in $\R^2$ and nondecreasing in each variable, then the
composite function $w=f(u_1,u_2)$ is ${\mathcal H}$--convex.
\end{lemma}
\begin{proof}
Assume first that $f\in C^2(\R^2)$, and set $ X_1=X,X_2=Y$. We
have
\[
X_jw=\sum_{p=1}^2\frac{\p f}{\p u_p}X_ju_p,
\]
\[
X_iX_jw=\sum_{p=1}^2\left(\frac{\p f}{\p u_p}X_iX_ju_p
+\sum_{q=1}^2\frac{\p^2 f}{\p u_q \p u_p}X_iu_qX_ju_p\right),
\]
and for every $h=(h_1,h_2)\in \R^2$
\[
\begin{split}
\langle \mathcal H(w) h,h\rangle &= \sum_{i,j=1}^2X_iX_jw
\,h_i\,h_j
\\
&=\sum_{p=1}^2\frac{\p f}{\p u_p}\langle \mathcal H(u_p)
h,h\rangle +\sum_{p,q=1}^2\frac{\p^2 f}{\p u_q \p
u_p}(\sum_{i=1}^2X_iu_qh_i)(\sum_{j=1}^2X_ju_ph_j)\\ &\geq 0,
\end{split}
\]
since $\mathcal H(u_p)$ is non negative definite and $\dfrac{\p
f}{\p u_p}\geq 0$ for $p=1,2$, and the matrix
\[
\left(\frac{\p^2 f}{\p u_q \p u_p}\right)_{p,q=1,2}
\]
is non negative definite.

If $f$ is only continuous, then given $h>0$ let
\[
f_h(x)=h^{-2}\int_{\R^2} \varphi\left(\frac{x-y}{h}\right)f(y)dy,
\]
where $\varphi\in C^\infty$ is nonnegative vanishing outside the
unit ball of $\R^2,$ and $\int \varphi =1.$ Since $f$ is convex,
then $f_h$ is convex and by the previous calculation
$w_h=f_h(u_1,u_2)$ is ${\mathcal H}$--convex. In particular, $w_h$
satisfies Proposition \ref{prop:H-conveximpliesconvexonlines} and
since $w_h\to w$ uniformly on compact sets as $h\rightarrow 0$, we
get that $w$ is ${\mathcal H}$--convex.
\end{proof}

\begin{proposition}\label{pro:osc}
Let $u\in C^2(\Omega)$ be ${\mathcal H}$--convex. For any compact
domain $\O'\Subset \O$ there exists a positive constant $C$
depending on $\Omega'$ and $\Omega$ and independent of $u$, such
that
\begin{equation}
\label{eq:osc} \int_{\O'}\{\det \mathcal H(u)+12 (u_t)^2\}\,dz\leq
C ({\rm osc}_\O u)^2.
\end{equation}
\end{proposition}
\begin{proof}
Given $\xi_0\in \O$ let $B_R=B_R(\xi_0)$ be a $d$--ball of radius
$R$ and center at $\xi_0$ such that $B_R\subset \O.$ Let $B_{\s
R}$ be the concentric ball of radius $\s R,$ with $0<\s<1.$
Without loss of generality we can assume $\xi_0=0,$ because the
vector fields $X$ and $Y$ are left invariant with respect to the
group of translations. Let $M=\max_{B_R} u$, then $u-M\leq 0$ in
$B_R$. Given $\varepsilon>0$ we shall work with the function
$u-M-\varepsilon <-\varepsilon$. In other words, by subtracting a
constant, we may assume $u<-\e$ in $B_R,$ for each given positive
constant $\e$; $\e$ will tend to zero at the end of the proof.

Define
\[
m_0=\inf_{B_R} u,
\]
and
\[
\varv(\xi)=\frac{m_0}{(1-\s^4)R^4}(R^4-\|\xi\|^4).
\]
Obviously $\varv=0$ on $\p B_R$ and $\varv=m_0$ on $\partial B_{\s
R}.$ We claim that $\varv$ is $\mathcal H$--convex in $B_R$ and
$\varv\leq m_0$ in $B_{\s R}.$ Setting $r=\|\xi\|^4$,
$h(r)=\dfrac{m_0}{(1-\s^4)R^4}(R^4-r)$, and following the
calculations in the proof of Proposition
\ref{prop:Hofdistanceiszero} we get
\[
\det \mathcal
H(\varv)=144(x^2+y^2)^2\left(\frac{m_0}{(1-\s^4)R^4}\right)^2\geq
0,
\]
 and
 \[
X^2 h = Y^2 h=-12\,(x^2+y^2)\frac{m_0}{(1-\s^4)R^4}\geq 0,
 \]
because $m_0$ is negative. Hence $\varv$ is $\mathcal H$--convex
in $B_R.$ Since $\varv-m_0=0$ on $\partial B_{\s R}$, it follows
from Proposition \ref{core} that $\varv\leq m_0$ in $B_{\s R}$. In
particular, $\varv\leq u$ in $B_{\s R}.$

Let $\rho\in C_0^\infty(\R^2)$, radial with support in the
Euclidean unit ball, $\int_{\R^2}\rho(x)\,dx=1$, and let
\begin{equation}\label{eq:functionfh}
f_h(x_1,x_2)=h^{-2}\,
\int_{\R^2}\rho((x-y)/h)\,\max\{y_1,y_2\}\,dy_1dy_2.
\end{equation}
We have that
\begin{enumerate}
\item
If $x_1>x_2$, then there exists $h_0>0$ and a neighborhood $V$ of
$(x_1,x_2)$ such that $f_h(y_1,y_2)=y_1$ for all $(y_1,y_2)\in V$
and for all $h\leq h_0$.\footnote{If $x_1>x_2$, then there exists
a cube $Q$ centered at $(x_1,x_2)$ such that if $(z_1,z_2)\in Q$
then $z_1>z_2$. Hence $x_1-y_1>x_2-y_2$ for all $|(y_1,y_2)|<h$
with $h$ sufficiently small. Then
\begin{align*}
f_h(x_1,x_2) &= h^{-2}\,\int_{|y|<h} \rho(y/h) (x_1-y_1)\,dy_1dy_2
=
x_1 -h^{-2}\,\int_{|y|<h} \rho(y/h) y_1\,dy_1dy_2\\
&= x_1 - h \int_0^1 t^2 \rho(t)\int_{S^1} y_1\,d\sigma(y)\,dt=x_1.
\end{align*}}
\item
There exists a positive constant $\alpha$ such that $f_h(x,x)=x +
\alpha \, h$ for all $h>0$ and for all $x\in \R$.\footnote{We have
\begin{align*}
f_h(x,x) &=
h^{-2}\,\int_{|y|<h} \rho(y/h) \max\{x-y_1,x-y_2\}\,dy_1dy_2\\
&= h^{-2}\,\int_{|y|<h} \rho(y/h)\left(x+
\max\{-y_1,-y_2\}\right)\,dy_1dy_2\\
&= x_1 + h^{-2}\,\int_{|y|<h} \rho(y/h)
\max\{-y_1,-y_2\}\,dy_1dy_2\\
&= x_1 + h^{-2}\,\int_{|y|<h} \rho(y/h)
\max\{y_1,y_2\}\,dy_1dy_2 \\
&= x_1 + h\,\int_{|y|<1} \rho(y)
\max\{y_1,y_2\}\,dy_1dy_2\\
&= x_1 + h\,\int_0^1 t^2  \rho(t)\int_{S^1}
\max\{y_1,y_2\}\,d\sigma(y)\,dt\\
&= x_1 + h\,\int_0^1 t^2  \rho(t)\int_{S^1}
\dfrac{|y_1-y_2|+ y_1 + y_2}{2}\,d\sigma(y)\,dt\\
&= x_1 + h\,\int_0^1 t^2  \rho(t)\int_{S^1}
\dfrac{|y_1-y_2|}{2}\,d\sigma(y)\,dt=x_1 + \alpha \,h.
\end{align*}}
\item
For all $h>0$, $f_h(\cdot,x_2)$ is nondecreasing for each $x_2$
and $f_h(x_1,\cdot)$ is nondecreasing for each $x_1$.
\end{enumerate}
Define
\[
w_h=f_h(u,\varv).
\]
From Lemma \ref{lem:comp} $w_h$ is $\mathcal H$--convex in $B_R$.
If $y\in B_{\sigma R}$ then $\varv(y)\leq u(y)$. If
$\varv(y)<u(y)$ then $f_h(u,\varv)(y)=u(y)$ for $h$ sufficiently
small; and if $\varv(y)=u(y)$, then $f_h(u,\varv)(y)= u(y)+\alpha
\,h$. Hence
\begin{align}\label{eq:uinsigmarcontrolled}
\int_{B_{\s R}} \{\det \mathcal H(u)+12(\p_tu)^2\}\,dz &=
\int_{B_{\s R}}
\{\det \mathcal H(w_h)+12((w_h)_t)^2\}\,dz\notag\\
&\leq \int_{B_R} \{\det \mathcal H(w_h)+12((w_h)_t)^2\}\,dz.
\end{align}
Now notice that $f_h(u,\varv)\geq \varv$ in $B_R$ for all $h$
sufficiently small. In addition, $u<0$ and $\varv=0$ on $\partial
B_R$ so $f_h(u,\varv)=0$ on $\partial B_R$. Then we can apply
Theorem \ref{thm:comparisonprinciple} to $w_h$ and $\varv$ to get
\[
\begin{split}
\int_{B_R} \{\det \mathcal H(w_h)+12(\p_tw_h)^2\}\,dz &\leq
\int_{B_R}
\{\det \mathcal H(\varv)+12(\varv_t)^2\}\,dz\\
&= 48\left( \frac{m_0}{(1-\s)
R^4}\right)^2\int_{B_R}(3(x^2+y^2)^2+t^2)\,dz\\
& = 48\left( \frac{m_0}{(1-\s)
}\right)^2\int_{B_1}(3(x^2+y^2)^2+t^2)\,dz.
\\
\end{split}
\]
This inequality combined with \eqref{eq:uinsigmarcontrolled}
yields
\[
\int_{B_{\s R}} \{\det \mathcal H(u)+12(\p_tu)^2\}\,dz\leq
C\,({m_0})^2 \leq C\, ({\rm osc}_{B_R} u + \varepsilon)^2.
\]
The inequality \eqref{eq:osc} then follows letting $\e\rightarrow
0$ and covering $\O'$ with balls.
\end{proof}
\begin{corollary}\label{corol}
Let $u\in C^2(\Omega)$ be ${\mathcal H}$--convex. For any compact
domain $\O'\Subset \O$ there exists a positive constant $C,$
independent of $u,$ such that
\begin{equation}
\label{eq:osc2} \int_{\O'}\det \mathcal H(u)\,dz\leq C ({\rm
osc}_\O u)^2.
\end{equation}
\end{corollary}
\begin{corollary}\label{corol2}
Let $u\in C^2(\Omega)$ be ${\mathcal H}$--convex. For any compact
domain $\O'\Subset  \O$ there exists a positive constant $C,$
independent of $u,$ such that
\begin{equation}
\label{eq:osc3} \int_{\O'}\trace {\mathcal H}_2(u)\,dz\leq C R^2
{\rm osc}_\O u.
\end{equation}
\end{corollary}

\subsection{Measure generated by an ${\mathcal H}$--convex function}
We shall prove that the notion $\int \det \mathcal H(u)+u_t^2$ can
be extended for continuous and ${\mathcal H}$--convex functions as
a Borel measure. We call this measure the ${\mathcal H}$--measure
associated with $u$, and we shall show that the map $u\in
C(\O)\rightarrow \mu(u)$ is weakly continuous on $C(\O).$
\begin{theorem}\label{measure}
Given an ${\mathcal H}$--convex function $u\in C(\O)$ there exists
a unique Borel measure $\mu(u)$ such that when $u\in C^2(\Omega)$,
\begin{equation}\label{eq:meas}
\mu(u)(E)=\int_E \{\det \mathcal H(u)+12u_t^2\}\,dz
\end{equation}
for any Borel set $E\subset \O.$ Moreover, if $u_k\in C(\Omega)$
are ${\mathcal H}$--convex, and $u_k\to u$ on compact subsets of
$\O,$ then $\mu(u_k)$ converges weakly to $\mu(u),$ that is,
\begin{equation}\label{eq:weakly}
\int_\O f\,d\mu(u_k) \rightarrow \int_\O f\,d\mu(u),
\end{equation}
for any $f\in C(\O)$ with compact support in $\O$.
\end{theorem}
\begin{proof} Let $u\in C(\O)$ be ${\mathcal H}$--convex, and let
$\{u_k\}\subset C^2(\O)$ be a sequence of ${\mathcal H}$--convex
functions converging to $u$ uniformly on compacts of $\O$. By
Proposition \ref{pro:osc}
\[
\int_{\O'}\{\det \mathcal H(u_k)+12 (\p_t u_k)^2\}\,dz
\]
are uniformly bounded, for every $\O'\Subset\O$, and hence a
subsequence of $(\det \mathcal H(u_k)+12 (\p_t u_k)^2)$ converges
weakly in the sense of measures to a Borel measure $\mu(u)$ on
$\O.$ We now prove that the map $u\in C(\Omega) \rightarrow
\mu(u)\in M(\O),$ the space of finite Borel measures on $\O$, is
well defined. Accordingly, let $\{\varv_k\}\subset C^2(\O)$ be
another sequence of ${\mathcal H}$--convex functions converging to
$u$ uniformly on compacts of $\O$. Assume $(\det \mathcal
H(u_k)+12 (\p_t u_k)^2)$ and $(\det \mathcal H(\varv_k)+12 (\p_t
\varv_k)^2)$ converge weakly to Borel measures $\mu, \mu'$
respectively. Let $B=B_R\Subset \O,$ and fix $\s \in (0,1).$ Let
$\eta\in C^2(\bar \O)$ be an ${\mathcal H}$--convex function such
that $\eta=0$ in $B_{\s R}$ and $\eta=1$ on $\p B_R$.\footnote{In
the $d$--ball $B_R(0)$, the function $\eta$ can be constructed as
follows. If $\varv(\xi)=\dfrac{1}{1-\sigma^4}
\left(\dfrac{\|\xi\|^4}{R^4}-\sigma^4 \right)$ and $f_h$ is the
function given by \eqref{eq:functionfh}, then define
$\eta(\xi)=f_h(\varv,0)$ with $h$ sufficiently small.} From the
uniform convergence of $\{u_k\}$ and $\{\varv_k\}$ towards $u$,
given $\e>0$ there exists $ k_\e\in \N$ such that
\[
-\frac{\e}{2}\leq u_k(x)-\varv_k(x) \leq \frac{\e}{2}, \quad
\text{for all $x\in \bar B$ and $k\geq k_\e$.}
\]
Hence
\[
u_k + \frac{\e}{2} \leq \varv_k+\e \eta
\]
on $\p B_R$ for $k\geq k_\e$. Define $\O_k=\{ \xi \in B_R:
u_k+\frac{\e}{2}>\varv_k+\e \eta\}.$ From Theorem
\ref{thm:comparisonprinciple} we have
\begin{align}\label{eqm1}
\int_{\O_k}\{\det \mathcal H(u_k)+12 (\p_t u_k)^2\}\,dz &\leq
\int_{\O_k}\det \mathcal H(\varv_k+\e \eta)+12 (\p_t \varv_k+\e
\p_t\eta)^2\notag\\
&\leq  \int_{B_R}\det \mathcal H(\varv_k)+12 (\p_t \varv_k)^2+\e^2 \, C \notag \\
&\qquad +\e \, C\, \int_{B_R}\left(\trace {\mathcal
H}_2(\varv_k)+|\p_t
\varv_k|\right) \notag\\
&\leq  \int_{B_R}\det \mathcal H(\varv_k)+12 (\p_t \varv_k)^2+\e^2\, C \notag\\
&\qquad +\e \,C\, \int_{B_R}\left(\trace {\mathcal
H}_2(\varv_k)+|\p_t \varv_k|^2+1\right)
\end{align}
and by Proposition \ref{pro:osc} and Corollary \ref{corol2} the
right hand side is bounded by
\[
\int_{B_R}\det \mathcal H(\varv_k)+12 (\p_t \varv_k)^2+\e \, C.
\]
By definition of $\Omega_k$ and since $\eta=0$ in $B_{\sigma R}$,
it follows that $B_{\s R}\subset \O_k$ and so by \eqref{eqm1} we
get
\begin{equation}\label{eqm2}
\int_{B_{\s R}}\det \mathcal H(u_k)+12 (\p_t u_k)^2\leq
\int_{B_R}\det \mathcal H(\varv_k)+12 (\p_t \varv_k)^2+\e \, C,
\end{equation}
and letting $k\rightarrow \infty,$ we get $\mu(B_{\s R})\leq
\mu'(B_R)+C\,\e$. Hence if $\e\rightarrow 0$ and $\s\rightarrow 1$
we obtain
\[
\mu(B)\leq \mu'(B).
\]
By interchanging $\{u_k\}$ and $\{\varv_k\}$ we get $\mu=\mu'.$

To prove \eqref{eq:weakly}, we first claim that it holds when
$u_k\in C^2(\Omega)$. Indeed, let $u_{k_m}$ be an arbitrary
subsequence of $u_k$, so $u_{k_m}\to u$ locally uniformly as $m\to
\infty$. By definition of $\mu(u)$, there is a subsequence
$u_{k_{m_j}}$ such that $\mu\left( u_{k_{m_j}} \right)\to \mu(u)$
weakly as $j\to \infty$. Therefore, given $f\in C_0(\Omega)$, the
sequence $\int_\Omega f\,d\mu(u_k)$ and an arbitrary subsequence
$\int_\Omega f\,d\mu(u_{k_m})$, there exists a subsequence
$\int_\Omega f\,d\mu(u_{k_{m_j}})$ converging to $\int_\Omega
f\,d\mu(u)$ as $j\to \infty$ and \eqref{eq:weakly} follows. For
the general case, given $k$ there exists $u_j^k\in C^2(\Omega)$
such that $u_j^k\to u_k$ locally uniformly as $j\to \infty$. By
definition of $\mu(u_k)$, there exists a subsequence $u_{j_m}^k$
such that $\mu\left( u_{j_m}^k \right)\to \mu(u_k)$ weakly as
$m\to \infty$. Let $f\in C_0(\Omega)$, $\supp f=K\subset
\Omega'\Subset \Omega$. There exists $m_1<m_2<\cdots $ such that
\[
|u_{j_{m_k}}^k(z)-u_k(z)|<1/k ,\qquad \text{for all $z\in
\Omega'$,}
\]
and
\[
\left| \int_\Omega f\, d\mu \left(u_{j_{m_k}}^k \right) -
\int_\Omega f\, d\mu(u_k) \right|<1/k,
\]
for $k=1,2,\cdots $. Hence $\varv_k=u_{j_{m_k}}^k\to u$ uniformly
in $\Omega'$ as $k\to \infty$, and so from the previous claim
\[
\int_\Omega f\,d\mu(\varv_k)\to \int_\Omega f\,d\mu(u), \quad
\text{as $k\to \infty$.}
\]
Therefore,
\begin{align*}
\left| \int_\Omega f\,d\mu(u_k) - \int_\Omega f\,d\mu(u)\right|
&\leq \left| \int_\Omega f\,d\mu(u_k) - \int_\Omega
f\,d\mu(\varv_k)\right| + \left| \int_\Omega f\,d\mu(\varv_k) -
\int_\Omega f\,d\mu(u)\right|
\\
&\leq \frac{1}{k}+ \left| \int_\Omega f\,d\mu(\varv_k) -
\int_\Omega f\,d\mu(u)\right|\to 0, \text{as $k\to \infty$},
\end{align*}
and the proof of the theorem is complete.
\end{proof}

\begin{corollary}\label{cor:generalcomparison}
If $u,\varv\in C(\bar \Omega)$ are $\mathcal H$--convex in
$\Omega$, $u=\varv$ on $\partial \Omega$ and $u\geq \varv$ in
$\Omega$, then $\mu(u)(\Omega)\leq \mu(\varv)(\Omega)$.

\end{corollary}

\subsection{Comparison principle for $\mathcal H$--measures}

\begin{theorem}\label{thm:comparisonpcpleformeasures}
Let $\Omega\subset \R^3$ be an open bounded set. If $u,\varv\in
C(\bar \Omega)$ are $\mathcal H$--convex in $\Omega$, $u\leq
\varv$ on $\partial \Omega$ and $\mu(u)(E)\geq \mu(\varv)(E)$ for
each $E\subset \Omega$ Borel set, then $u\leq \varv$ in $\Omega$.
\end{theorem}

\begin{proof}
Assume $0\in \Omega$, $\Delta=\diam (\Omega)$, $\e>0$, and
$u_\e(x,y,t)=u(x,y,t)+\e\,(x^2+y^2 -\Delta^2)$. We have $x^2+y^2
-\Delta^2<0$ for $(x,y,t)\in \bar \Omega$, so $u_\e < u\leq v$ in
$\partial \Omega$. Suppose there exists $(x_o,y_o,t_o)\in \Omega$
such that $u(x_o,y_o,t_o)>\varv(x_o,y_o,t_o)$. Hence the set
$D=\{(x,y,t)\in \Omega:u_\e(x,y,t)>\varv(x,y,t)\}$ is non empty
for all $\e$ sufficiently small. In addition, $\bar D\cap \partial
\Omega=\emptyset$. So $\bar D\subset \Omega$ and $u_\e=\varv$ on
$\partial D$. By Corollary \ref{cor:generalcomparison} we get
$\mu(u_\e)(D)\leq \mu(\varv)(D)$. On the other hand, there exist
$u_k\in C^2(\Omega)$ $\mathcal H$--convex in $\Omega$ such that
$u_k\to u$ uniformly on compact subsets of $\Omega$. Let
$u_{k,\e}(x,y,t)= u_k(x,y,t) + \e \, (x^2+y^2 -\Delta^2)$. We have
from \eqref{eq:detofuepsilon} that
\begin{align*}
\int_D \{\det \mathcal H(u_{k,\e}) + (u_{k,\e})_t^2\}\,dz &=
\int_D \{\det \mathcal H(u_k)+2\e\, \trace \mathcal H(u_k)
+ 4\e^2 + (u_k)_t^2\}\,dz\\
&\geq \mu(u_k)(D) + 4\e^2 \, |D|.
\end{align*}
Letting $k\to \infty$ we get from Theorem \ref{measure} that
$\mu(u_\e)(D)\geq  \mu(u)(D) + 4\e^2 \, |D|>\mu(u)(D)$ obtaining a
contradiction.
\end{proof}

\bibliographystyle{alpha}

\end{document}